# Low Rank Approximation in $G_0W_0$ Calculations

Meiyue Shao[*]    Lin Lin[†]    Chao Yang[‡]    Fang Liu[§]

Felipe H. da Jornada[¶]    Jack Deslippe[‖]    Steven G. Louie[**]

May 7, 2016

## Abstract

The single particle energies obtained in a Kohn–Sham density functional theory (DFT) calculation are generally known to be poor approximations to electron excitation energies that are measured in transport, tunneling and spectroscopic experiments such as photo-emission spectroscopy. The correction to these energies can be obtained from the poles of a single particle Green's function derived from a many-body perturbation theory. From a computational perspective, the accuracy and efficiency of such an approach depends on how a self energy term that properly accounts for dynamic screening of electrons is approximated. The $G_0W_0$ approximation is a widely used technique in which the self energy is expressed as the convolution of a non-interacting Green's function ($G_0$) and a screened Coulomb interaction ($W_0$) in the frequency domain. The computational cost associated with such a convolution is high due

[*]Computational Research Division, Lawrence Berkeley National Laboratory, Berkeley, CA 94720, USA (`myshao@lbl.gov`)

[†]Department of Mathematics, University of California, Berkeley and Computational Research Division, Lawrence Berkeley National Laboratory, Berkeley, CA 94720, USA (`linlin@math.berkeley.edu`).

[‡]Computational Research Division, Lawrence Berkeley National Laboratory, Berkeley, CA 94720, USA (`cyang@lbl.gov`).

[§]School of Statistics and Mathematics, Central University of Finance and Economics, Beijing 100081, China (`fliu@cufe.edu.cn`).

[¶]Department of Physics, University of California at Berkeley & Materials Science Division, Lawrence Berkeley National Laboratory, Berkeley, CA 94720, USA (`jornada@berkeley.edu`)

[‖]NERSC, Lawrence Berkeley National Laboratory, Berkeley, CA 94720, USA (`jdeslippe@lbl.gov`).

[**]Department of Physics, University of California at Berkeley & Materials Science Division, Lawrence Berkeley National Laboratory, Berkeley, CA 94720, USA (`sglouie@berkeley.edu`)

to the high complexity of evaluating $W_0$ at multiple frequencies. In this paper, we discuss how the cost of $G_0W_0$ calculation can be reduced by constructing a low rank approximation to the frequency dependent part of $W_0$. In particular, we examine the effect of such a low rank approximation on the accuracy of the $G_0W_0$ approximation. We also discuss how the numerical convolution of $G_0$ and $W_0$ can be evaluated efficiently and accurately by using a contour deformation technique with an appropriate choice of the contour.

# 1 Introduction

Electron excitations in molecules and solids that, for example, can be measured in photo-emission spectroscopy experiments can be modeled by the theory of single particle Green's functions $G(\boldsymbol{r}, \boldsymbol{r}'; t)$. In the frequency domain, the poles of $G$, which correspond to the excitation energies of a molecule or solid, are eigenvalues of a single particle Hamiltonian $H$ that contains a correction to the ground-state Hamiltonian

$$H_{\text{KS}} = -\nabla^2 + V_{ion} + V_{Hartree} + V_{xc}$$

obtained in Kohn–Sham density functional theory, where $V_{ion}$ is the electron-ion interaction potential, $V_{Hartree}$ is the Hartree potential, and $V_{xc}$ is the exchange–correlation potential, respectively. The correction is frequency (denoted by $\omega$) dependent. To be specific, we can write $H$ as

$$H(\omega) = H_{\text{KS}} - V_{xc} + \Sigma(\omega),$$

where $\Sigma$ is an energy dependent self energy operator that accounts for dynamic screening of electrons. The operator $H_{\text{KS}}$ is a self-adjoint operator and all its eigenvalues are real. However, in the presence of the frequency-dependent self energy operator $\Sigma(\omega)$, $H(\omega)$ is in general not a self-adjoint operator.

If $(\varepsilon_j^{\text{KS}}, \psi_j^{\text{KS}})$ is a normalized eigenpair of $H_{\text{KS}}$, it follows from the first order eigenvalue perturbation theory that the corresponding eigenvalue of $H$ can be approximated by

$$\varepsilon_j \approx \varepsilon_j^{\text{KS}} + \langle \psi_j^{\text{KS}} | \left[ \Sigma(\varepsilon_j) - V_{xc} \right] | \psi_j^{\text{KS}} \rangle, \tag{1}$$

where $\langle \cdot | \cdot | \cdot \rangle$ is the Dirac bra–ket notation. These approximate eigenvalues can be obtained by solving the nonlinear scalar equation

$$\omega = \varepsilon_j^{\text{KS}} + \langle \psi_j^{\text{KS}} | \left[ \Sigma(\omega) - V_{xc} \right] | \psi_j^{\text{KS}} \rangle$$

by an appropriate numerical method. These methods would require evaluating $\Sigma(\omega)$ at multiple frequencies.

One of the challenges in this type of excited states calculation is to obtain a computable approximation to the self energy term $\Sigma(\varepsilon)$ in (1), which does not have a closed form expression. We will use $\Sigma(\omega)$ and $\Sigma(\varepsilon)$ interchangeably below since energy and frequency are closely related quantities. In this paper, we assume all Kohn–Sham orbitals $\psi_j^{\text{KS}}(\boldsymbol{r})$ are real. This corresponds to the case for molecules and for solids with inversion symmetry. To simplify our discussion, we do not distinguish an integral operator and its kernel unless otherwise clarified. Throughout this paper, we use $\eta$ to denote an infinitesimal real positive number. The final quasi-particle energy obtained by solving Eq. (1) should not depend on the value of $\eta$, but the sign of $\eta$ is important due to the analytic structure of various operators in the frequency domain in the complex plane.

A widely used approximation, called the $G_0W_0$ *approximation* [1, 4, 5, 9, 16, 18, 29], expresses the self energy operator $\Sigma(\boldsymbol{r}, \boldsymbol{r}'; \omega)$ as the convolution of a non-interacting Green's function $G_0$ and a screened Coulomb interaction $W_0$:

$$\Sigma(\boldsymbol{r}, \boldsymbol{r}'; \omega) = \frac{\mathrm{i}}{2\pi} \int_{\mathbb{R}} G_0(\boldsymbol{r}, \boldsymbol{r}'; \omega + \omega') W_0(\boldsymbol{r}', \boldsymbol{r}; \omega') \mathrm{e}^{\mathrm{i}\omega'\eta} \, \mathrm{d}\omega'. \tag{2}$$

Here the term $\mathrm{e}^{\mathrm{i}\omega'\eta}$ with infinitesimally small $\eta > 0$ is needed due to causality. We refer readers to a recent work on the mathematical analysis of the GW method [6] for more details.

The time ordered Green's function is defined as

$$G_0(\boldsymbol{r}, \boldsymbol{r}'; \omega) = \sum_{j=1}^{n_v} \frac{\psi_j^{\text{KS}}(\boldsymbol{r}) \psi_j^{\text{KS}}(\boldsymbol{r}')}{\omega - \varepsilon_j^{\text{KS}} - \mathrm{i}\eta} + \sum_{j=n_v+1}^{n} \frac{\psi_j^{\text{KS}}(\boldsymbol{r}) \psi_j^{\text{KS}}(\boldsymbol{r}')}{\omega - \varepsilon_j^{\text{KS}} + \mathrm{i}\eta}, \tag{3}$$

where $n_v$ is the number of occupied (or valence) states, $n$ is the total number of single particle states of a discretized $H_{\text{KS}}$. Here the signs of $\mathrm{i}\eta$ in Eq. (3) have physical meanings, and reflect the different propagation directions in time for electrons and holes (see, e.g., [1, 6]).

The screened Coulomb interaction takes the form

$$W_0(\boldsymbol{r}, \boldsymbol{r}'; \omega) = \int_{\mathbb{R}^3} \epsilon^{-1}(\boldsymbol{r}, \boldsymbol{r}''; \omega) v(\boldsymbol{r}'', \boldsymbol{r}') \, \mathrm{d}\boldsymbol{r}'', \tag{4}$$

where $v(\boldsymbol{r}, \boldsymbol{r}')$ is the bare Coulomb interaction, and $\epsilon(\boldsymbol{r}, \boldsymbol{r}'; \omega)$ is the dielectric function defined by

$$\epsilon(\boldsymbol{r}, \boldsymbol{r}'; \omega) = \delta(\boldsymbol{r}, \boldsymbol{r}') - \int_{\mathbb{R}^3} v(\boldsymbol{r}, \boldsymbol{r}'') \chi_0(\boldsymbol{r}'', \boldsymbol{r}'; \omega) \, \mathrm{d}\boldsymbol{r}'', \tag{5}$$

and $\chi_0(\boldsymbol{r}, \boldsymbol{r}'; \omega)$ is the irreducible polarizability function

$$\chi_0(\boldsymbol{r}, \boldsymbol{r}'; \omega) = 2 \sum_{i=1}^{n_v} \sum_{j=n_v+1}^{n} \psi_i^{\text{KS}}(\boldsymbol{r}) \psi_j^{\text{KS}}(\boldsymbol{r}) \psi_i^{\text{KS}}(\boldsymbol{r}') \psi_j^{\text{KS}}(\boldsymbol{r}') \\ \times \left( \frac{1}{\omega - (\varepsilon_j^{\text{KS}} - \varepsilon_i^{\text{KS}}) + i\eta} - \frac{1}{\omega + (\varepsilon_j^{\text{KS}} - \varepsilon_i^{\text{KS}}) - i\eta} \right). \quad (6)$$

A direct $G_0W_0$ calculation by evaluating the self energy operator $\Sigma(\omega)$ at multiple frequencies is computationally expensive. The high cost mainly comes from two sources: 1) Eq. (2) involves integrating a function in $\omega'$ over the entire real axis. The function to be integrated is singular at multiple frequencies $\omega'$ in the limit of $\eta = 0$. The singularity of $G_0$ can be seen directly from Eq. (3), which has $n$ poles that are $\eta$-distance away from the real axis. It can also be shown that the screened Coulomb interaction $W_0$ has many poles close to the real axis. Hence, directly replacing the integral by a numerical quadrature, even when performed carefully, requires evaluating the integrand at many frequencies. This is not only expensive, but can also be potentially unstable numerically when additional approximations of $W_0$ are made to reduce the computational cost. 2) In order to compute $W_0$, we need to perform a double summation over products of occupied and unoccupied wavefunctions $\psi_i^{\text{KS}}(\boldsymbol{r}) \psi_j^{\text{KS}}(\boldsymbol{r})$ to obtain $\chi_0$ first. For large systems, the complexity of $\chi_0$ calculation is $\mathcal{O}(n^4)$. Furthermore, an explicit calculation of $\chi_0$ requires unoccupied Kohn–Sham wavefunctions $\psi_j^{\text{KS}}$, $j = n_v + 1, \ldots, n$ to be computed in addition to the occupied states returned from a ground-state calculation.

In order to reduce the computational cost, efforts have been made to replace $\chi_0$ by a low rank approximation, using, for example, a truncated singular value decomposition (SVD). However, because the singular values of $\chi_0$ decay slowly, a relatively large number of singular vectors need to be computed to achieve sufficient accuracy. The cost of computing these singular vectors is not significantly lower than the cost of forming $\chi_0$ directly.

It has been observed in [15, 22] that the singular values of $v\chi_0$ decay more rapidly. In particular, the largest few singular values, can be orders of magnitude larger than other singular values. This is mainly due to contribution from bare Coulomb interaction $v$ at small wave numbers. Therefore, instead of constructing a low rank approximation to $\chi_0$, one can construct a low rank approximation to $v\chi_0$. Such a low rank approximation can be otained efficiently by using an iterative method that only requires multiplying $\chi_0$ with a number of vectors. It is well known that this type of multiplication can be carried out without constructing $\chi_0$ explicitly first [13, 25, 34]. Because $\chi_0$ is a linear response operator, the product $\chi_0$ and a vector can be obtained via density functional perturbation theory [3], which requires solving a number of linear systems of equations with a shifted Kohn–Sham Hamiltonian $H_{\text{KS}}$ as the coefficient matrix. These linear systems of equations can in turn be solved iteratively by using, for example, preconditioned Krylov subspace methods.

A low rank approximation to $v\chi_0$ yields a low rank approximation to the frequency dependent part of $W_0$, which we will denote by $W_p$. Although this low rank approximation is not an optimal approximation to $W_p$ in the matrix 2-norm, the error introduced by such an approximation can be bounded. An optimal low rank approximation to $W_p$ should be obtained by performing a truncated SVD of $W_p$ directly, which is more costly. However, we should note that the integrand to be computed in (2) is a pointwise multiplication (Hadamard product) of $G_0$ and $W_0$. Therefore, an optimal low rank approximation to $W_p$ does not necessarily give the optimal approximation to the integrand. Other low approximation schemes, such as the one based on setting long wavelength matrix elements of the Fourier space representation of $v\chi_0$ to zero, which is used in many existing packages such as BerkeleyGW [9], are also valid.

We should point out that neither $v\chi_0$ nor $W_p$ is low rank in the traditional sense, i.e., the smallest singular values of these operators are not tiny compared to the largest singular values. The reason that a low rank approximation of $W_p$ can be used in $G_0W_0$ calculation is due to the fact that the quantity we calculate is ultimately a first order correction to Kohn–Sham eigenvalues. Thus, approximations to $v\chi_0$ and $W_p$ are acceptable as long as the error introduced by the approximation is comparable or smaller than the error associated with the first order perturbation theory.

In this paper, we will examine different ways to construct a low rank approximation to $W_p$ and analyze the effect of such low rank approximation on the accuracy of the self energy calculation. Numerical experiments are performed to study the effectiveness of using a low rank approximation of $W_p$ in the $G_0W_0$ calculation. To assess the accuracy related exclusively to the low rank approximation of $W_p$, we assume all linear systems that must be solved to construct a low rank approximation of $W_p$ are solved accurately, and operations that involve working with $G_0$ are performed accurately also. Therefore, the only source of error in the evaluation of the integrand in (2) originates from the low rank approximation. We also ensure that numerical integration is performed accurately in (2) using a special technique called contour deformation. This technique is discussed in detail in Section 3.

To simplify notation, we assume all quantities have been discretized in real space. The operators or functions $G_0$, $W_0$, $v$, $\chi_0$, and $\Sigma$ will be treated as matrices. For instance, the integral (4) will simply be written as a matrix–matrix multiplications $W_0 = \epsilon^{-1}v$. The spatial variables $\boldsymbol{r}$ and $\boldsymbol{r}'$ are often omitted unless otherwise clarified.

Our convention also assumes that a function of a single variable becomes a column vector after it is discretized.

The rest of this paper is organized as follows. In Section 2 we briefly review the technique of splitting the self energy into frequency independent and frequency dependent parts. In Section 3 we propose a new contour deformation strategy by introducing the concept of residue free frequency. A numerical integration scheme using Gaussian quadrature rules is also discussed. In Section 4 we present several low rank approximations to $W_p$, and provide error estimates. Finally we show by a numerical example the effectiveness of proposed strategies.

## 2 The separation of exchange from correlation in the self energy

From the definition of the dielectric function in Eq. (5), i.e., $\epsilon(\omega) = I - v\chi_0(\omega)$, one can derive by the Sherman–Morrison–Woodbury (SMW) formula (see, for instance, [8]) that the inverse of the dielectric function has the form

$$\epsilon^{-1}(\omega) = I + v\chi(\omega), \tag{7}$$

where $\chi$ is the reducible polarizability that satisfies the equation

$$\chi(\omega) = \chi_0(\omega) + \chi_0(\omega)v\chi(\omega),$$

or, equivalently,

$$\chi = [I - \chi_0(\omega)v]^{-1}\chi_0(\omega). \tag{8}$$

When all Kohn–Sham orbitals $\psi_j$'s are real, both the Green's function $G_0$ and $W_0$ are *complex symmetric* for $\omega \in \mathbb{C}$, i.e.,

$$G_0(\boldsymbol{r}, \boldsymbol{r}'; \omega) = G_0(\boldsymbol{r}', \boldsymbol{r}; \omega), \quad W_0(\boldsymbol{r}, \boldsymbol{r}'; \omega) = W_0(\boldsymbol{r}', \boldsymbol{r}; \omega).$$

Note that due to the infinitesimal complex frequency $i\eta$, neither $G_0$ nor $W_0$ are Hermitian operators. The symmetry of $W_0$ follows from

$$W_0(\omega) = \epsilon^{-1}(\omega)v = v + v\chi(\omega)v, \tag{9}$$

where $\chi(\omega)$ can be shown to be complex symmetric due to Eq. (8) and the fact that $\chi_0(\omega)$ is complex symmetric.

In order to simplify the convolution in Eq. (2), it is often convenient to use (9) to split $W_0$ into the frequency independent bare Coulomb interaction and the frequency dependent interaction $v\chi(\omega)v$, which we denote by $W_p(\omega)$. While the bare Coulomb interaction is constant in the frequency domain, $W_p(\omega)$ decays to 0 as $|\omega| \to \infty$. Hence the $G_0W_0$ approximation of $\Sigma$ becomes

$$\Sigma(\omega) = \frac{\mathrm{i}}{2\pi} \int_{\mathbb{R}} G_0(\omega + \omega') \odot v \mathrm{e}^{\mathrm{i}\omega'\eta} \, \mathrm{d}\omega' + \frac{\mathrm{i}}{2\pi} \int_{\mathbb{R}} G_0(\omega + \omega') \odot W_p(\omega') \mathrm{e}^{\mathrm{i}\omega'\eta} \, \mathrm{d}\omega', \quad (10)$$

where $\odot$ is used to denote an element-wise multiplication (i.e., Hadamard product). Here we have used the complex symmetric property of $W_0$.

By using the residue theorem and the expression (3) for the Green's function, we can rewrite the first term in (10) as

$$\Sigma_X \equiv -\left[\sum_{j=1}^{n_v} \psi_j^{\mathrm{KS}} \left(\psi_j^{\mathrm{KS}}\right)^{\mathsf{T}}\right] \odot v.$$

It follows from the trace identity

$$a^{\mathsf{T}}(A \odot B)b = \mathrm{trace}\left(\mathrm{Diag}(a) \cdot A \cdot \mathrm{Diag}(b) \cdot B^{\mathsf{T}}\right) \quad (11)$$

for vectors $a$, $b$ and matrices $A$, $B$, where $\mathrm{Diag}(a)$ denotes a diagonal matrix with the vector $a$ on its diagonal, that we can express the $\Sigma_X$ contribution to $\langle\psi_i^{\mathrm{KS}}|\Sigma|\psi_i^{\mathrm{KS}}\rangle$ as

$$\langle\psi_i^{\mathrm{KS}}|\Sigma_X|\psi_i^{\mathrm{KS}}\rangle = -\mathrm{trace}\left[\sum_{j=1}^{n_v} \rho_{ij}^{\mathrm{KS}} \left(\rho_{ij}^{\mathrm{KS}}\right)^{\mathsf{T}} v\right] = -\sum_{j=1}^{n_v} \left(\rho_{ij}^{\mathrm{KS}}\right)^{\mathsf{T}} v \rho_{ij}^{\mathrm{KS}}, \quad (12)$$

where $\rho_{ij}^{\mathrm{KS}} = \psi_i^{\mathrm{KS}} \odot \psi_j^{\mathrm{KS}}$ and $v^{\mathsf{T}} = v$. Note that $\rho_{ij}^{\mathrm{KS}}$ is a column vector here since $\psi_i^{\mathrm{KS}}$ denotes the discretized $\psi_i^{\mathrm{KS}}(\boldsymbol{r})$ in the real space and is a column vector. Because the type of calculation involved in (12), which is frequency independent, is similar to that required in a Hartree–Fock exchange energy calculation, (12) is often referred to as the exchange part of the self energy. This is what the subscript $X$ in $\Sigma_X$ stands for. The remaining contribution, which is defined by the second term of (10), is referred to as the correlation part of the self energy and denoted by $\Sigma_C$. In [21], the numerical integration required to evaluate $\Sigma_C$ is carried out carefully on the real axis directly. We will discuss an alternative way to perform such a numerical integration when a low rank representation of $W_p$ is available in Section 3.

Numerical results indicate that $\Sigma_X$ is often larger compared to $\Sigma_C$. However, it is the relatively small contribution from $\Sigma_C$ (usually on the order of eV) that gives the quantitative prediction power of the $G_0W_0$ method for a variety of electron excitation properties.

**Algorithm 1:** Compute the correlation energy correction $\langle \psi_i^{\text{KS}} | \Sigma_C(\omega) | \psi_i^{\text{KS}} \rangle$ associated with $\varepsilon_i^{\text{KS}}$ for a given $\omega$.

---

**Input:** Ground-state Kohn–Sham eigenpairs $(\varepsilon_j^{\text{KS}}, \psi_j^{\text{KS}})$, for $j = 1, 2, \ldots, n$; the frequency $\omega$ at which the correlation energy is to be evaluated; the quasiparticle energy index $i$ (e.g., HOMO or LUMO)

**Output:** $\langle \psi_i^{\text{KS}} | \Sigma_C(\omega) | \psi_i^{\text{KS}} \rangle$

1: **for** each $\omega' \in$ {frequencies used to perform the integration} **do**
2:     Construct $\chi_0(\omega')$;
3:     Compute $\epsilon(\omega')$;
4:     Evaluate $W_p(\omega') = \epsilon^{-1}(\omega')v - v$;
5: **end for**
6: Perform numerical integration over $\omega'$ to obtain an approximate $\Sigma_C(\omega)$;
7: Compute $\langle \psi_i^{\text{KS}} | \Sigma_C(\omega) | \psi_i^{\text{KS}} \rangle$.

---

A brief outline of the major steps required to compute $\langle \psi_i^{\text{KS}} | \Sigma_C(\omega) | \psi_i^{\text{KS}} \rangle$ for a specific $i$ is given in Algorithm 1. Because the integrand in $\Sigma_C$ contains sharp peaks and troughs on the real axis, many $\omega'$ values need to be generated in order to achieve sufficient accuracy when we perform numerical integration. For each $\omega'$, the computational complexity is dominated by Step 2 of the algorithm, which is $\mathcal{O}(n^2 n_v n_c)$, where $n_c \equiv n - n_v$ is the number of unoccupied Kohn–Sham states.

## 3 Numerical integration via contour deformation

The difficulty of performing numerical integration of $\langle \psi_i^{\text{KS}} | \Sigma_C(\omega) | \psi_i^{\text{KS}} \rangle$ on the real axis can be overcome by using a technique that reformulates the improper integration as part of a contour integral along a closed contour $\Omega$ constructed on the complex plane. By using the Cauchy integral theorem, we can express the integral to be evaluated in terms of an integral along $\Omega \setminus (-\infty, +\infty)$ if the contour $\Omega$ does not enclose any poles of $G_0$ or $W_p$. This technique is often referred to as the *contour deformation* technique [14, 20].

## 3.1 The choice of contour

It follows from (9) that $W_0(\omega')$ and $W_p(\omega')$ have identical poles. When $\eta \neq 0$, it is not difficult to verify that all poles of $W_p(\omega')$ are in the second and fourth quadrants of the complex plane, and are symmetric with respect to the origin. If there exists a $\delta_W$ such that

$$\delta_W := \min_{\omega'}\{|\operatorname{Re}\omega'|\colon \omega' \text{ is a pole of } W_p\} > 0, \tag{13}$$

that is, if the distance between the poles of $W_p$ and the imaginary axis is nonzero, there is a stripe of nonzero width in which $W_p(\omega')$ is analytic. In this case, it is possible to choose a contour for $\omega'$ in the complex plane as shown in Figure 1 along which $G_0(\omega + \omega') \odot W_p(\omega')$ integrates to zero because the contour does not contain any poles of $G_0$ or $W_p$. Because $|G_0(\omega')|$ and $|W_p(\omega')|$ behave like $|\omega'|^{-1}$ and $|\omega'|^{-2}$ respectively in the limit of $|\omega'| \to \infty$ [10], the contribution to the contour integral from the two arcs (cyan color in Figure 1) vanish as the radius of the arcs increases towards $\infty$. Therefore, integrating $G_0(\omega+\omega') \odot W_p(\omega')$ along the real axis (red color in Figure 1) is equivalent to integrating $-G_0(\omega + \omega') \odot W_p(\omega')$ along the vertical path parallel to the imaginary axis (magenta color in Figure 1). This path consists of frequencies $\mu - \omega + i\zeta$, where $\zeta \in (-\infty, +\infty)$. Because the poles of $G_0$ and $W_p$ are sufficiently far from such a path, the integrand $G_0(\omega + \omega') \odot W_p(\omega')$ is a well behaved function along the path. As a result, it can be integrated accurately by an appropriate numerical quadrature scheme with a relatively small number of quadrature points.

We should point out that the assumption given by (13) is justified for molecules with a finite HOMO–LUMO gap

$$\varepsilon_g := \varepsilon_{n_v+1} - \varepsilon_{n_v}$$

centered near zero. Although the poles of $W_p$ are not exactly the same as those of $\chi_0$, it can be shown that they are close to the poles of $\chi_0$. As a result, $\delta_W > 0$ is likely to hold for molecules. Our numerical results for the $SiH_4$ molecule to be presented later indicate that $\delta_W$ can even be larger than $\varepsilon_g$. Hence, as a first approximation, one can use $\varepsilon_g$ to estimate the value of $\delta_W$.

We should also emphasize that the existence of a contour that does not enclose any poles of $G_0(\omega + \omega')$ or $W_p(\omega')$ depends on the choice of $\omega$ as well as the positions of the poles of $G_0(\omega + \omega')$ and $W_p(\omega')$. When such a contour exists, we refer to the corresponding $\omega$ a *residue free frequency* in the $G_0W_0$ theory.

We now characterize the range of residue free frequencies $\omega$ along the real axis. We denote the highest occupied state $\varepsilon_{n_v}^{\text{KS}}$ by $\varepsilon_{\text{HOMO}}^{\text{KS}}$ and the lowest unoccupied state

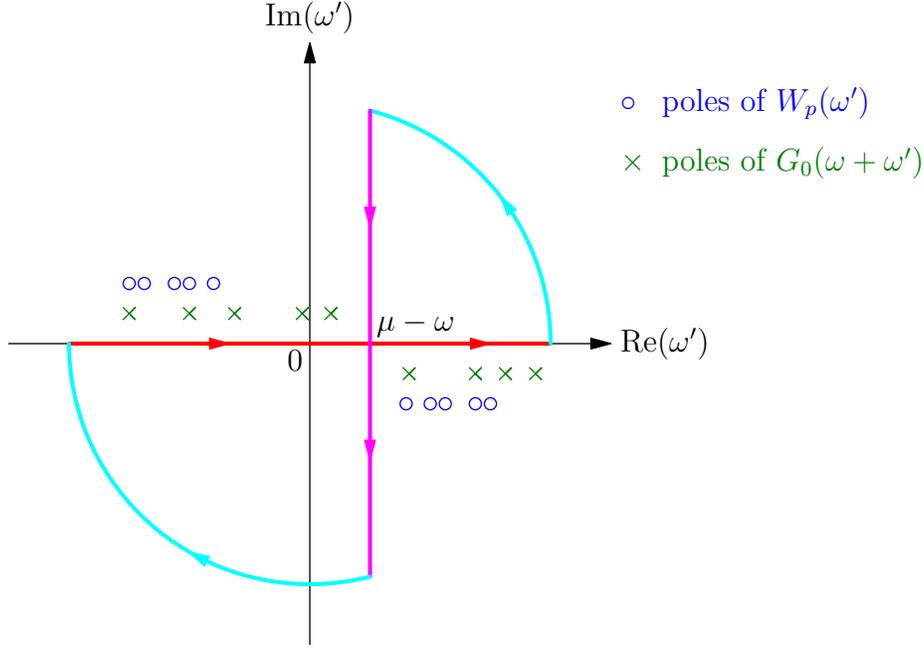

Figure 1: A contour that bypasses all poles of $G_0(\omega + \omega')v\chi(\omega')v$.

$\varepsilon^{\text{KS}}_{n_v+1}$ by $\varepsilon^{\text{KS}}_{\text{LUMO}}$. Define

$$\omega'_{\text{lb}} = \max\left\{-\delta_W, \varepsilon^{\text{KS}}_{\text{HOMO}} - \omega\right\}, \quad \omega'_{\text{ub}} = \min\left\{\varepsilon^{\text{KS}}_{\text{LUMO}} - \omega, \delta_W\right\}.$$

If $\omega'_{\text{lb}} < \omega'_{\text{ub}}$, we may choose the vertical path of integration $\omega' = (\omega'_{\text{lb}} + \omega'_{\text{ub}})/2 + i\zeta$, where $\zeta \in (-\infty, +\infty)$ as the vertical path of integration. Equivalently, the condition $\omega'_{\text{lb}} < \omega'_{\text{ub}}$ can be rewritten as

$$\varepsilon^{\text{KS}}_{\text{HOMO}} - \delta_W < \omega < \varepsilon^{\text{KS}}_{\text{LUMO}} + \delta_W. \tag{14}$$

When the condition (14) does not hold (Figure 2 shows such an example), it is not possible to construct a simple contour like the one shown in Figure 1 to exclude all poles of $G_0$ and $W_p$. Hence the condition (14) defines the set of residue free frequencies.

We remark that if $\omega$ is not a residue free frequency, it does not mean that the $G_0W_0$ calculation cannot proceed. If we can find a path parallel to the imaginary axis that does not contain any poles of $G_0$ or $W_p$, the integration along such a path can be corrected by adding residues associated with the poles enclosed by the contour. For example, if the poles of $G_0(\omega + \omega')$ and $W_p(\omega')$ are distributed as shown in Figure 2,

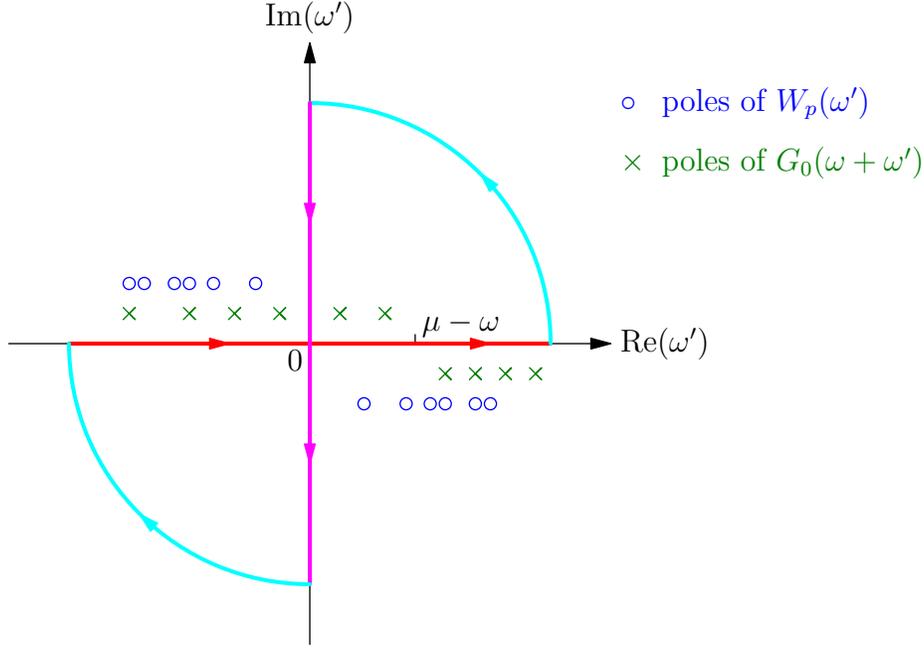

Figure 2: An example that (14) does not hold.

we can evaluate $\Sigma_C(\omega)$ as

$$\Sigma_C(\omega) = -\sum_{j=1}^{p} \lim_{z \to z_j} G_0(z-z_j) \odot W_p(z)(z-z_j) - \frac{1}{2\pi} \int_{-\infty}^{+\infty} G_0(\omega+\mathrm{i}\zeta) \odot W_p(\mathrm{i}\zeta) \,\mathrm{d}\zeta, \quad (15)$$

where $z_j$, $j = 1, 2, \ldots, p$ are the poles of $G_0(\omega + \omega')$ that are enclosed by the contour shown in the figure.

The evaluation of the first term in (15) involves computing $W_p(z_j)$ at a few poles of $G_0(\omega + \omega')$, that is, at $z_j = \varepsilon_{n_v-j+1}^{\mathrm{KS}} - \omega + \mathrm{i}\eta$ that are inside the contour. The integral in the second term of (15), which is performed along a path parallel to the imaginary axis, can be evaluated numerically by using an appropriate quadrature rule. In fact, the contour deformation technique presented in existing literature almost always chooses the imaginary axis [4, 14, 15, 20] as the vertial path of integration. For small molecules, we expect such an approach to work well, since the poles of $G_0(\omega + \omega')$ are not very clustered. However, for systems with an increasing number of atoms, the poles of $G_0$ can become clustered. In the case of solids, densely clustered eigenvalues form continuous energy bands. In such a case, it becomes increasingly likely that some poles of $G_0$ can be very close to the imaginary axis on which numerical integration is performed. It can also be the case that the residue term to be evaluated has a pole close to a pole of $W_p$. In both cases, the evaluation of $G_0(\omega + \omega')W_0(\omega')$ can become difficult, especially when the spectral function of $G_0$ is not available and when $W_p$ is approximated by a low rank operator.

## 3.2 Gaussian quadrature

Regardless how we choose the vertical path of the contour, we need to integrate (22) along such a path numerically. The numerical integration scheme we use depends on how the integrand (22) behaves on such a path.

Denote by $\bar{z}$ the complex conjugation of a complex number $z$, and $\bar{A}$ the elementwise complex conjugation of a matrix $A$, then from the definition of $G_0$ and $W_p$, when $\operatorname{Im}\omega' \neq 0$, we have

$$G_0(\boldsymbol{r},\boldsymbol{r}';\overline{\omega'}) = \overline{G_0(\boldsymbol{r},\boldsymbol{r}';\omega')}, \qquad W_p(\boldsymbol{r},\boldsymbol{r}';\overline{\omega'}) = \overline{W_p(\boldsymbol{r},\boldsymbol{r}';\omega')}. \qquad (16)$$

This symmetry property suggests that the integration along $\omega' = \omega_s + \mathrm{i}\zeta$, where $\zeta \in (-\infty, +\infty)$, can be reduced to integrating along the path that is in the upper half of the complex plane because

$$\begin{aligned}
\Sigma_C(\omega) &= -\frac{1}{2\pi} \int_{-\infty}^{+\infty} G_0(\omega + \omega_s + \mathrm{i}\zeta) \odot W_p(\omega_s + \mathrm{i}\zeta) \,\mathrm{d}\zeta \\
&= -\frac{1}{\pi} \int_0^{+\infty} G_0(\omega + \omega_s + \mathrm{i}\zeta) \odot W_p(\omega_s + \mathrm{i}\zeta) \,\mathrm{d}\zeta.
\end{aligned} \qquad (17)$$

Through the change of variable

$$\zeta = \frac{\xi}{1-\xi}, \qquad (18)$$

we turn the integration from 0 to $+\infty$ to that from 0 to 1 (see Figure 3), i.e.,

$$\Sigma_C(\omega) \xrightarrow{\zeta=\xi/(1-\xi)} -\frac{1}{\pi} \int_0^1 G_0(\omega + \omega_s + \mathrm{i}\zeta) \odot W_p(\omega_s + \mathrm{i}\zeta) \frac{1}{(1-\xi)^2} \,\mathrm{d}\xi. \qquad (19)$$

Due to the asymptotic decay $|G_0(\omega')| \sim |\omega'|^{-1}$ and $|W_p(\omega')| \sim |\omega'|^{-2}$ as $|\omega'| \to \infty$, it can be shown that the integrand of Eq. (19) has finite and nonzero limit both at $\xi = 0$ and $\xi = 1$. This observation suggests that (19) can be approximated by a Legendre–Gauss–Radau (LGR) quadrature. Because the integrand is well behaved, it is likely that only a small number of quadrature points is needed along the integration path. Our numerical results shown in Section 5 confirm this prediction. However, when the path of integration is close to poles of $G_0$ or $W_p$, the number of quadrature points may increase significantly in order to achieve the desired accuracy.

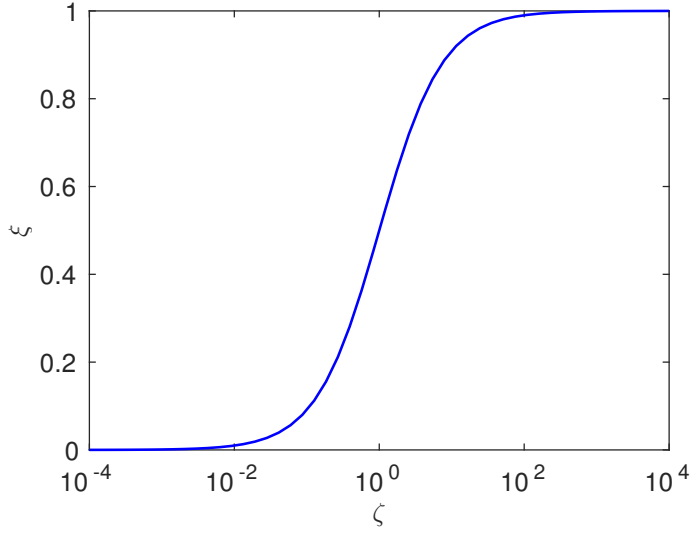

Figure 3: The curve of the transformation (18).

## 3.3 Estimating the poles of $W_p$

Choosing an appropriate contour to perform numerical integration of (22) requires estimating the position of the pole of $W_p(\omega')$ with the smallest positive real part. We define $\Phi$ to be the matrix that contains products of $\psi_i^{\text{KS}}$ and $\psi_j^{\text{KS}}$, previously denoted by $\rho_{ij}^{\text{KS}}$, as its columns for $i = 1, 2, \ldots, n_v$ and $j = n_v + 1, \ldots, n$, that is, $\Phi = \left(\rho_{ij}^{\text{KS}}, \ldots\right)$, and $D$ be the diagonal matrix that contains $\varepsilon_j^{\text{KS}} - \varepsilon_i^{\text{KS}}$ on it diagonal. It is well known [7, 33] that the poles of $W_p$ are the eigenvalues of the matrix

$$H_c = \begin{pmatrix} D + K & K \\ -K & -(D + K) \end{pmatrix}, \tag{20}$$

where $K = \Phi^T v \Phi$.

The dimension of this matrix is $2n_v n_c \times 2n_v n_c$. It can be shown that the desire eigenvalue of (20) can be obtained by computing the smallest eigenvalue of the matrix $D^2 + 2D^{1/2}KD^{1/2}$. An iterative method such as the locally optimal block preconditioned conjugate gradient (LOBPCG) method [19] can be used to accomplish such a task. However, this calculation does require the use of unoccupied states $\psi_j^{\text{KS}}$, $j = n_v + 1, \ldots, n$. When these unoccupied states are not available or difficult to compute, we may use $\varepsilon_g$ as a rough estimate of $\delta_W$ since the poles of $W_p$ and $\chi_0$ are usually close to each other.

# 4 Low rank approximation of $W_p$

Because (1) holds only approximately, it may not be necessary to evaluate $\langle \psi_i^{\text{KS}}|\Sigma(\omega)|\psi_i^{\text{KS}}\rangle$ to full accuracy. In fact, as long as the error in the evaluation of $\langle \psi_i^{\text{KS}}|\Sigma(\omega)|\psi_i^{\text{KS}}\rangle$ is on the order $\mathcal{O}(\|V_{xc} - \Sigma(\omega)\|^2)$, the computed solution to the approximate equation (1) may be acceptable.

## 4.1 Low rank approximation via truncated SVD

We are interested in approximations to (10) that can lower the cost of $\langle \psi_i^{\text{KS}}|\Sigma(\omega)|\psi_i^{\text{KS}}\rangle$ computation. One way to reduce such cost is to replace $W_p(\omega')$ in (10) by a low rank approximation obtained from, for example, a truncated singular value decomposition (SVD) of the form

$$W_p(\omega') \approx U(\omega')S(\omega')V(\omega')^*, \tag{21}$$

where $S$ is a $k \times k$ diagonal matrix that contains $k$ largest singular values of $W_p$ for some $k \ll n$, and $U$ and $V$ are $n \times k$ matrices containing the corresponding right and left singular vectors. When $\omega'$ is real and the $\eta \to 0^+$ limit is taken, $W_p(\omega')$ is in fact Hermitian, and the SVD can be replaced by an eigenvalue decomposition. However, we will use SVD to account for the general case in which $\omega'$ is a complex number. The choice of complex $\omega'$ is required when contour deformation is used to evaluate the integral in (10).

Using (11) and (21), we can rewrite $\langle \psi_i^{\text{KS}}|\Sigma_C(\omega)|\psi_i^{\text{KS}}\rangle$ as

$$\langle \psi_i^{\text{KS}}|\Sigma_C(\omega)|\psi_i^{\text{KS}}\rangle \approx \frac{\mathrm{i}}{2\pi}\int_{\mathbb{R}} \text{trace}\left[\text{Diag}(\psi_i^{\text{KS}})G_0(\omega+\omega')\text{Diag}(\psi_i^{\text{KS}})U(\omega')S(\omega')V(\omega')^*\right]\mathrm{d}\omega'. \tag{22}$$

If we define $U_{\psi_i^{\text{KS}}}(\omega') \equiv \text{Diag}(\psi_i^{\text{KS}})U(\omega')$ and $V_{\psi_i^{\text{KS}}}(\omega') = \text{Diag}(\psi_i^{\text{KS}})V(\omega')$, the expression (22) can be simplified to yield

$$\langle \psi_i^{\text{KS}}|\Sigma_C(\omega)|\psi_i^{\text{KS}}\rangle \approx \frac{\mathrm{i}}{2\pi}\int_{\mathbb{R}} \text{trace}\left[V_{\psi_i^{\text{KS}}}(\omega')^* G_0(\omega+\omega') U_{\psi_i^{\text{KS}}}(\omega') S(\omega')\right]\mathrm{d}\omega'. \tag{23}$$

The integrand in (23) can be evaluated as the trace of a $k \times k$ matrix, which can be calculated efficiently as long as

1. the truncated SVD of $W_p$ can be easily obtained;

2. $G_0(\omega+\omega')U(\omega')$ can be computed efficiently.

A truncated SVD of $W_p$ can be computed iteratively by using a subspace iteration (see, for example, [27]) to extract the dominant singular values and the corresponding right and left singular vectors of $W_p$. This procedure involves repeatedly multiplying $W_p$ and $W_p^*$ with a block of vectors denoted by $X \in \mathbb{C}^{n \times k}$. If we use the expression $W_p(\omega) = [I - v\chi_0(\omega)]^{-1} v - v$ in $W_p(\omega)X$, these multiplications would require solving linear systems of the form

$$[I - v\chi_0(\omega)]Y_1 = X_1, \quad \text{and} \quad [I - v\chi_0(\omega)]^* Y_2 = X_2. \tag{24}$$

These linear systems of equations can be solved iteratively by using methods such as (block) GMRES [28], QMR [12], or BiCGSTAB [35]. Each step of these methods requires multiplying $v\chi_0(\omega)$ and/or its conjugate transpose with a block of vectors.

It has been recognized [2, 15] that one of the attractive features of using iterative methods to solve (24) is that multiplying $\chi_0(\omega)$ with a vector $g$ does not require $\chi_0(\omega)$ to be constructed explicitly in advance. It is well known that $\chi_0(\omega)g$ can be evaluated as

$$\chi_0 g = \sum_{i=1}^{n_v} \text{Diag}(\psi_i^{\text{KS}}) \Delta \psi_i,$$

where $\Delta \psi_i$ is the solution to the equation

$$\left[ H_{\text{KS}} - (\varepsilon_i^{\text{KS}} - \omega)I \right] \Delta \psi_i = -P_v^\perp \text{Diag}(\psi_i^{\text{KS}})g, \tag{25}$$

for $i = 1, 2, \ldots, n_v$. Here $P_v^\perp = I - \sum_{i=1}^{n_v} \psi_i^{\text{KS}}(\psi_i^{\text{KS}})^\mathsf{T}$ is the orthogonal projection operator that projects against the invariant subspace spanned by $\{\psi_i\}_{i=1}^{n_v}$. Eq. (25) is often referred to as the *Sternheimer equation*. In addition to avoiding constructing $\chi_0(\omega)$ explicitly, which requires $\mathcal{O}(n^2 n_v n_c)$ operations, we do not need to compute any unoccupied (empty) Kohn–Sham orbitals. The solution of the Sternheimer equation is the most time consuming step. One possible way to solve Eq. (25) for multiple frequency is by the Lanczos method (e.g., see [15]). In this work we solve Eq. (25) by direct inversion. This is inefficient for large systems, and we will study the efficient numerical method more systematically in a future publication.

Once $U(\omega')$ becomes available, $U_{\psi_i^{\text{KS}}}(\omega')$ can be easily computed by row scaling. The use of unoccupied states can also be avoided in the evaluation of $G_0(\omega + \omega')U_{\psi_i^{\text{KS}}}(\omega')$ if we treat $G_0(\omega + \omega')$ simply as the inverse of $(\omega + \omega')I - H_{\text{KS}}$ instead of working with its spectral decomposition. In the limit of $\eta \to 0^+$, the product $G_0(\omega + \omega')U_{\psi_i^{\text{KS}}}(\omega')$ can be computed by solving a linear system of the form

$$[(\omega + \omega')I - H_{\text{KS}}] X = U_{\psi_i^{\text{KS}}}(\omega') \tag{26}$$

**Algorithm 2:** Compute $\langle \psi_i^{\text{KS}} | \Sigma_C(\omega) | \psi_i^{\text{KS}} \rangle$ using a low rank approximation of $W_p$

**Input:** Ground-state Kohn–Sham eigenpairs $(\varepsilon_j^{\text{KS}}, \psi_j^{\text{KS}})$, for $j = 1, 2, \ldots, n_v$; the frequency $\omega$ at which the correlation energy is to be evaluated; the quasiparticle energy index $i$ (e.g. HOMO or LUMO)

**Output:** $\langle \psi_i^{\text{KS}} | \Sigma_C(\omega) | \psi_i^{\text{KS}} \rangle$

1: **for** each $\omega' \in \{$frequencies used to perform the integration of (23)$\}$ **do**
2:   Construct a low rank approximation to $W_p$, i.e., $W_p(\omega') = U(\omega')S(\omega')V^*(\omega')$, where $U, V$ are $n \times k$, $S$ is a $k \times k$ diagonal matrix for $k \ll n$;
3:   $U_{\psi_i^{\text{KS}}}(\omega') \leftarrow \text{Diag}(\psi_i^{\text{KS}})U$;
4:   $V_{\psi_i^{\text{KS}}}(\omega') \leftarrow \text{Diag}(\psi_i^{\text{KS}})V$;
5:   Compute $G_0(\omega + \omega')U_{\psi_i^{\text{KS}}}(\omega')$ by solving a number of linear systems of the form (26);
6:   Take the trace of $V_{\psi_i^{\text{KS}}}(\omega')G_0(\omega + \omega')U_{\psi_i^{\text{KS}}}(\omega')S(\omega')$;
7: **end for**
8: Perform numerical integration over $\omega'$ to obtain an approximation to $\langle \psi_i^{\text{KS}} | \Sigma_C(\omega) | \psi_i^{\text{KS}} \rangle$;

iteratively by using MINRES/SYMMLQ [23] or COCG [30].

When a nonzero $\eta$ is used in (3), we can use the alternative expression of the Green's function

$$G_0(\omega) = P_v \left[(\omega - i\eta)I - H_{\text{KS}}\right]^{-1} P_v + P_v^\perp \left[(\omega + i\eta)I - H_{\text{KS}}\right]^{-1} P_v^\perp,$$

and solve two separate linear systems of the form

$$P_v \left[(\omega + \omega' - i\eta)I - H_{\text{KS}}\right] P_v X_1 = P_v U_{\psi_i^{\text{KS}}}(\omega'), \tag{27}$$

$$P_v^\perp \left[(\omega + \omega' + i\eta)I - H_{\text{KS}}\right] P_v^\perp X_2 = P_v^\perp U_{\psi_i^{\text{KS}}}(\omega'). \tag{28}$$

The sum of $X_1$ and $X_2$ gives the product $G_0(\omega + \omega')U_{\psi_i^{\text{KS}}}(\omega')$.

The main steps of using a low rank approximation of $W_p$ to obtain approximate $\langle \psi_i^{\text{KS}} | \Sigma_C(\omega) | \psi_i^{\text{KS}} \rangle$ for a particular frequency $\omega$ are outlined in Algorithm 2. The overall complexity of this procedure is dominated by the complexity of Steps 2 of the algorithm which requires solving at least $kn_v$ linear systems of the form (25).

The right-hand side of each of these equation takes at least $\mathcal{O}(nn_v^2)$ operations to construct. The multiplication of $H_{\text{KS}}$ with a vector in each iteration generally has a lower complexity if implemented efficiently. Therefore, assuming the number of iterations required to solve each linear system is a relatively small constant, the overall complexity of this approach is $\mathcal{O}(kn_v^3 n)$, which is not necessarily lower than that associated with Algorithm 1. However, when $n_v$ and $k$ are significantly smaller than $n$, this approach is more efficient. Furthermore, this approach is more amenable to parallel implementation. We will defer a more detailed discussion on how to use an iterative solver to obtain a truncated SVD of $W_p$ to a future publication.

## 4.2 Alternative low rank approximations of $W_p$

As we can see from the discussion in the previous subsection, constructing a truncated SVD of $W_p$ requires solving a set of linear systems of equations in a nested fashion, i.e., the solution of (24) requires the solution of (25). As a result, the error introduced in iterative solutions of these linear systems may accumulate and make the approximation less effective unless each linear system is solved accurately, which can be costly.

On the other hand, even though a truncated SVD of $W_p$ provides the best approximation to $W_p$ in 2-norm, it may not necessarily be the best low rank decomposition for approximating the integrand in (10) because $W_p$ appears as one of the two factors of a Hadamard product.

There are other alternative low rank approximation strategies that are simpler, less costly and potentially more effective. For instance, we can perform a truncated SVD on $v\chi(\omega) = \epsilon^{-1}(\omega) - I$ instead to obtain $W_p \approx \hat{U}\hat{S}(v\hat{V})^*$, where

$$v\chi(\omega) \approx \hat{U}\hat{S}\hat{V}^*.$$

Note that, in this low rank approximation, we retain one of the bare Coulomb factors $v$ in $W_p$, but relax the orthogonality of the $v\hat{V}$ factor. However, this approach would still require us to solve equations of the form (24).

We can also construct a truncated SVD approximation to $v\chi_0$, i.e., $v\chi_0(\omega) \approx \hat{U}\hat{S}\hat{V}^*$, and then obtain a low rank approximation of $v\chi$ via the Sherman–Morrison–Woodbury formula. In this case, the approximation to $W_p$ can be written as

$$W_p = \left[\epsilon^{-1}(\omega) - I\right]v \approx \hat{U}\hat{O}^{-1}(v\hat{V})^*, \tag{29}$$

where $\hat{O} = \hat{S}^{-1} - \hat{V}^*\hat{U}$. Note that $\hat{O}$ is an $k \times k$ matrix for $k \ll n$. Inverting such a matrix should not be too costly when $k$ is relatively small. In this approach, we only need to solve Sternheimer equations (25).

In principle, we need to construct different low rank approximations to $W_p(\omega)$ for different $\omega$. In the recent work by Govoni and Galli [15], it is suggested that the small dimensional subspaces constructed from low rank approximation of $W_p(0)$ can be used to construct low rank approximations of $W_p(\omega)$ for other nonzero $\omega$ values.

Another type of low rank approximation to $W_p(\omega)$ that is currently implemented in many software packages such as BerkeleyGW [9] is based on truncating the bare Coulomb operator $v$ in the reciprocal space. If $F$ is the matrix representation of a discrete Fourier transform, and $\hat{v}_t$ is the matrix representation of the bare Coulomb operator truncated in the reciprocal space, i.e., $\hat{v}_t$ is a $k \times k$ diagonal matrix obtained by extracting the $k$ largest diagonal elements of $F^*vF$, this low rank approximation has the form
$$W_p \approx (FE_t\hat{v}_t) \left[(\hat{\chi}_0)_t^{-1}(\omega) - \hat{v}_t\right]^{-1} (FE_t\hat{v}_t)^*, \tag{30}$$
where $E_t$ is an $n \times k$ matrix that contains $k$ selected columns of the $n \times n$ identity matrix, and $(\hat{\chi}_0)_t$ is the projection of $F^*\chi_0 F$ onto the subspace represented by $E_t$, i.e.,
$$(\hat{\chi}_0)_t = E_t^T F^* \chi_0(\omega) F E_t.$$
The selection of $k$ columns in $E_t$ is done conformally with the selection of the diagonal elements of $F^*vF$.

To see why (30) corresponds to a truncation of the bare Coulomb operation in the reciprocal space, let us write the reciprocal space representation of $v^{1/2}\chi_0(\omega)v^{1/2}$, which is similar to $v\chi_0(\omega)$, as
$$(F^*v^{1/2}F)F^*\chi_0(\omega)F(F^*v^{1/2}F) = \hat{v}^{1/2}\hat{\chi}_0(\omega)\hat{v}^{1/2},$$
where $\hat{v} \equiv F^*vF$ is a diagonal matrix. By setting diagonal elements of $\hat{v}^{1/2}$ that are smaller than a prescribed threshold to zeros, we effectively obtain a low rank approximation to $v^{1/2}\chi_0 v^{1/2}$ in the reciprocal space:
$$F^*v^{1/2}\chi_0(\omega)v^{1/2}F \approx E_t\hat{v}_t^{1/2}[\hat{\chi}_0]_t(\omega)\hat{v}_t^{1/2}E_t^T, \tag{31}$$
where $E_t$ and $\hat{v}_t$ are as defined previously.

Transforming (31) back to real space and multiplying it from the left and right by $v^{1/2}$ and $v^{-1/2}$ yields the approximation of $v\chi_0$:
$$v\chi_0 \approx \left(v^{1/2}FE_t\hat{v}_t^{1/2}\right)[\hat{\chi}_0]_t(\omega)\left(v^{-1/2}FE_t\hat{v}_t^{1/2}\right)^*.$$

It follows from the SMW formula, and the fact that $v^{1/2}FE_t\hat{v}_t^{1/2} = FE_t\hat{v}_t$, that $W_p$ can be approximated by (30).

Finally, we remark that in the GW community the resolution of identity (RI) technique is also frequently used to construct low-rank approximations. For discussions with the RI technique, we refer the readers to [11, 26, 36].

## 4.3 Error estimates

In the following we estimate the amount of error produced in the evaluation of

$$\langle \psi_i^{\text{KS}}|G_0(\omega+\omega') \odot W_p(\omega')|\psi_i^{\text{KS}}\rangle, \tag{32}$$

when a low rank approximation to $W_p(\omega')$ is used. For the sake of simplicity, we assume that all linear equations to be solved are solved exactly, and the only source of error originates in the low rank approximation to $W_p$.

If the rank-$k$ approximation of $W_p$ is obtained by the truncated SVD on $W_p$, the amount of error we introduced in $W_p$ can be quantified by a perturbation $\Delta W_p$ with $\|\Delta W_p\|_2 \leq \sigma_{k+1}(v\chi v)$, where $\sigma_i(\cdot)$ denotes the $i$th largest singular value. As a result, it follows from the inequality $\|A \odot B\| \leq \|A\|\|B\|$ (see, for example, [17, Section 5.1]), that

$$\left|\langle \psi_i^{\text{KS}}|G_0 \odot (W_p + \Delta W_p)|\psi_i^{\text{KS}}\rangle - \langle \psi_i^{\text{KS}}|G_0 \odot W_p|\psi_i^{\text{KS}}\rangle\right| \leq \|G_0\|\|\Delta W_p\| \leq \sigma_{k+1}(v\chi v)\|G_0\|. \tag{33}$$

Similarly, if the truncated SVD is applied to $v\chi$ instead of $v\chi v$, we have an error estimate

$$\left|\langle \psi_i^{\text{KS}}|G_0 \odot (W_p + \Delta W_p)|\psi_i^{\text{KS}}\rangle - \langle \psi_i^{\text{KS}}|G_0 \odot W_p|\psi_i^{\text{KS}}\rangle\right| \leq \sigma_{k+1}(v\chi)\|v\|\|G_0\|. \tag{34}$$

To estimate the error of truncating $v\chi_0$, we make use of the perturbation bound (see, for example, [8])

$$\left\|\epsilon^{-1} - (\epsilon+\Delta\epsilon)^{-1}\right\| \leq \frac{\|\epsilon^{-1}\|\|\epsilon^{-1}\Delta\epsilon\|}{1-\|\epsilon^{-1}\Delta\epsilon\|}.$$

When $\|\Delta\epsilon\| = \sigma_{k+1}(v\chi_0)$ is sufficiently small (say, $\|\epsilon^{-1}\Delta\epsilon\| \leq 1/2$), we have

$$\left\|\epsilon^{-1} - (\epsilon+\Delta\epsilon)^{-1}\right\| \leq 2\sigma_{k+1}(v\chi_0)\|\epsilon^{-1}\|^2$$

and hence

$$\left|\langle \psi_i^{\text{KS}}|G_0 \odot (W_p + \Delta W_p)|\psi_i^{\text{KS}}\rangle - \langle \psi_i^{\text{KS}}|G_0 \odot W_p|\psi_i^{\text{KS}}\rangle\right| \leq 2\sigma_{k+1}(v\chi_0)\|\epsilon^{-1}\|^2\|v\|\|G_0\|. \tag{35}$$

Let us denote the upper bounds in (33)–(35) by

$$\begin{aligned} E_1(k) &= \sigma_{k+1}(v\chi v)\|G_0\|, \\ E_2(k) &= \sigma_{k+1}(v\chi)\|v\|\|G_0\|, \\ E_3(k) &= 2\sigma_{k+1}(v\chi_0)\|\epsilon^{-1}\|^2\|v\|\|G_0\|, \end{aligned}$$

respectively. It is well known that truncated SVD provides the best rank-$k$ approximation by the Eckart–Young theorem [32]. Therefore we have

$$E_1(k) \leq E_2(k) \leq E_3(k).$$

However, we remark that the error bounds provided here are very pessimistic, and a higher error bound does not necessarily imply a larger error.

So far we have assumed that the integral in (23) can be evaluated efficiently and accurately by an appropriate numerical quadrature rule. We will see in the subsequent section that the integrand in (23) appears to be sufficiently smooth using the contour deformation technique. Thus we can use a Gaussian quadrature rule to evaluate the numerical integration accurately.

## 4.4 Methods for solving Sternheimer equations

The Sternheimer equations that appear in both (25) and the equations results from evaluating $G_0(\omega + \omega')U(\omega')$ can be solved iteratively. Because the coefficient matrix is frequency dependent, we need to solve one set of equations for each frequency. This can become expensive when the number of quadrature points used in the quadrature rule is large. However, if $\omega'$ is away from the real axis, the coefficient matrix is likely to be well conditioned. In this case, we can first generate a frequency independent subspace by using an $m$-step Lanczos algorithm for $m \ll n$ to yield

$$H_{\text{KS}}V_m = V_m T_k + \beta v_{m+1} e_m^\mathsf{T}, \qquad \text{with} \qquad (V_m, v_{m+1})^\mathsf{T}(V_m, v_{m+1}) = I,$$

where $T_m$ is symmetric and tridiagonal. Approximation to the solution of a linear system of the form

$$(H_{\text{KS}} - \omega I)x = b, \tag{36}$$

can be obtained by standard Galerkin techniques that require solving either

$$(T_m - \omega I)g = \|b\|e_1,$$

or

$$\min_g \left\| \begin{pmatrix} T_m \\ \beta e_m^\mathsf{T} \end{pmatrix} g - \|b\|e_1 \right\|,$$

and forming $x = V_m g$ at different frequencies. Existing implementations of the GMRES, MINRES, QMR, and BiCGSTAB algorithms can be easily modified to account for frequency dependency without applying the standard solver to each linear system associated with a different frequency.

Because we need to solve multiple linear systems with different right hand sides in (23) and (25), we can use recently developed Krylov subspace recycling techniques [24, 31] to further reduce the computational cost. In a parallel implementation, we can map different right hand sides to different processor groups to improve parallel scalability.

We should note that when $\omega'$ is close to a pole of $G_0(\omega + \omega')$ or $W_p(\omega')$, which can occur when a poor contour is chosen to perform the numerical integration (23), or when $\omega$ is not a residue free frequency and the vertical path of integration comes close to a pole of $G_0(\omega + \omega')$, the linear system (36) can become ill-conditioned. In this case, many iterations may be required to solve this linear system unless an effective preconditioner is identified and used. If we terminate the Lanczos procedure too early by keeping $m$ too small, large error may be introduced in the numerical integration.

The computational details of developing an efficient iterative solver for Sternheimer equations will be described in a separate paper.

## 5 Numerical examples

In this section, we demonstrate the effect of low rank approximation on the accuracy of the self energy calculation for a small molecule, silane ($SiH_4$). We use the MATLAB toolbox KSSOLV [37] to obtain ground-state Kohn–Sham eigenvalues and wavefunctions (orbitals). The number of valence (occupied) states is $n_v = 4$. The planewave (kinetic energy) cutoff we used in this calculation is Ecut = 15 Ry, resulting in a total of $n_g = 939$ planewaves employed to represent each Kohn–Sham orbital. The number of unoccupied states is thus $n_c = n_g - n_v = 935$.

### 5.1 The poles of $G_0$ and $W_p$

Because this problem is relatively small, we can compute the poles of $W_p(\omega')$ explicitly by solving an eigenvalue problem of dimension $n_v n_c = 3780$. The smallest positive pole is at $\omega' = \delta_W = 3.3952\,\text{eV}$. This is also the value of $\delta_W$ that appears in (14). For this problem, the gap between the HOMO and LUMO energies is $\varepsilon_g = \varepsilon_{\text{LUMO}}^{\text{KS}} - \varepsilon_{\text{HOMO}}^{\text{KS}} = 3.1948\,\text{eV}$. Thus it serves as a good lower bound for $\delta_W$. As we discussed in Section 3.1, for $\omega$'s that satisfy (14), there exists a path of integration parallel to the imaginary axis that is free of any poles of $G_0$ and $W_p$.

In Figure 4, we plot both the poles of $G_0(\omega + \omega')$ and $W_p(\omega')$ for two different $\omega$ values: $\omega = (\varepsilon_{n_v} + \varepsilon_{n_v+1})/2$ and $\omega = \varepsilon_{n_v+1} + \delta_W$. To see the relative positions of these poles clearly, we use $\eta_1 = 0.01\,\text{eV}$ for the poles of $G_0$ and $\eta_2 = 0.02\,\text{eV}$ for the poles of $W_p$. We observe that when $\omega = (\varepsilon_{n_v} + \varepsilon_{n_v+1})/2$, the imaginary axis is clearly free of any poles of $G_0$ and $W_p$. When $\omega \geq \varepsilon_{n_v+1} + \delta_W$, the poles of $G_0$ and $W_p$ start to overlap. Any contour of the form shown in Figures 1 or 2 will enclose some poles of either $G_0$ or $W_p$.

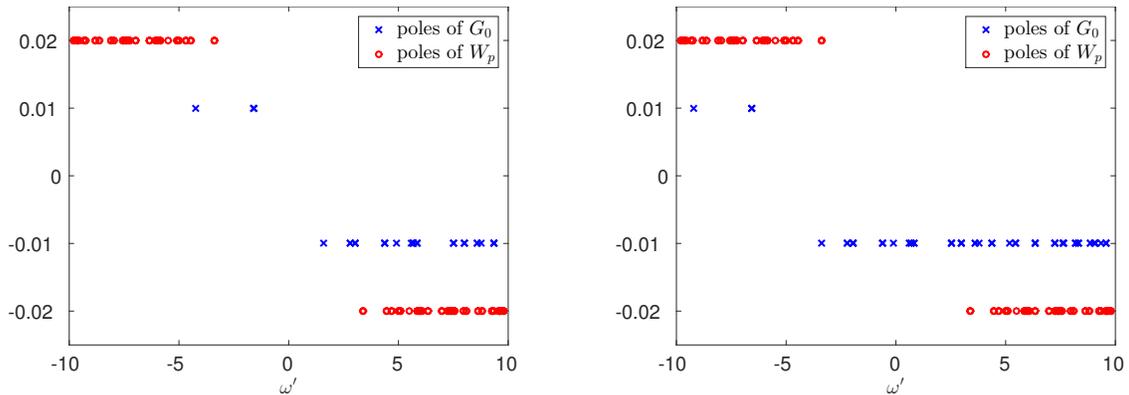

Figure 4: The poles of $G_0(\omega + \omega')$ and $W_p(\omega')$ for $\omega = (\varepsilon_{n_v} + \varepsilon_{n_v+1})/2$ (left) and $\omega = \varepsilon_{n_v+1} + \delta_W$ (right).

## 5.2 Contour choices

When $\omega = (\varepsilon_{n_v} + \varepsilon_{n_v+1})/2 = 0.2630 \, \text{eV}$, there is a reasonably wide vertical stripe, $-\varepsilon_g/2 < \text{Re}(\omega') < \varepsilon_g/2$, that is free of poles of $G_0$ and $W_p$. In particular, we can choose the imaginary axis as part of the contour of integration to evaluate (22).

Figure 5 shows how the real part of the integrand in (22) changes for $i = 4$ (HOMO), along the imaginary axis. Since the integrand is an even function as shown in (16), we plot only the integrand in the upper half of the imaginary axis. A change of variable of the form (18) is used to map $\zeta \in [0, +\infty)$ to $\xi \in [0, 1)$.

To illustrate the benefit of choosing a contour that is sufficiently far away from any poles of $G_0$ and $W_p$, we plot the real part of the integrand in (22) for $i = 4$ (HOMO) along the vertical paths, $\omega' = (1 - \delta)\varepsilon_g/2 + \mathrm{i}\zeta$, $\zeta \in [0, +\infty)$, of several different contours with $\delta \in \{10^0, 10^{-1}, 10^{-2}, 10^{-3}\}$ in Figure 6. Again, we use the change of variable (18) to map $\zeta \in [0, +\infty)$ to $\xi \in [0, 1)$.

The contour associated with $\delta = 1$ is the best one in the sense of being furthest away from any pole of $G_0$ and $W_p$. The other contours become progressively closer to a pole of $G_0$, as $\delta$ becomes smaller.

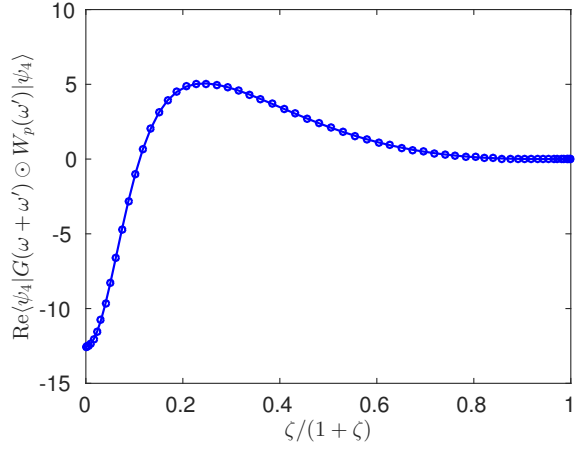

Figure 5: The real part of $\langle\psi_4|G_0(\omega+\omega')\odot W_p(\omega')|\psi_4\rangle$ along the imaginary axis. A change of variable of variable of the form (18) is used to map $[0,+\infty)$ to $[0,1)$.

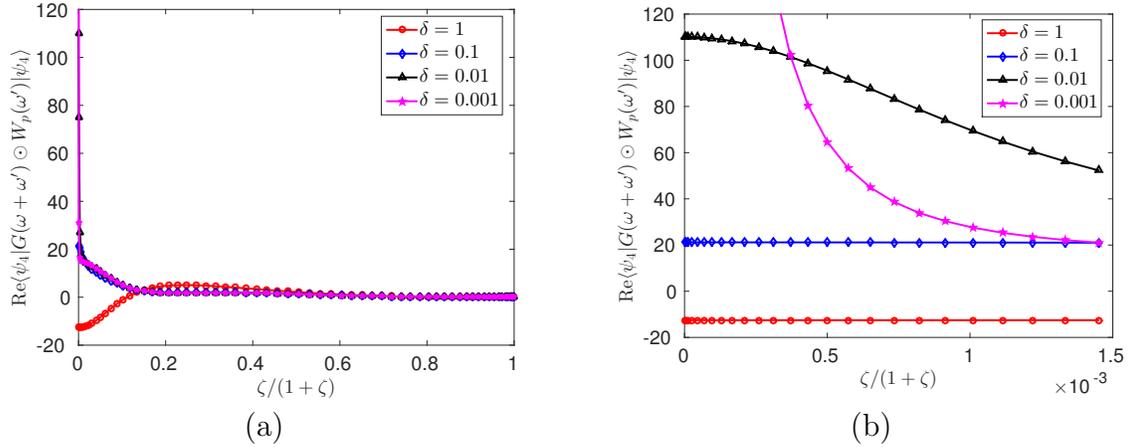

(a)          (b)

Figure 6: (a) The real part of $\langle\psi_4|G(\omega+\omega')\odot W_p(\omega')|\psi_4\rangle$ along the lines $\omega' = (1-\delta)\varepsilon_g/2 + \mathrm{i}\zeta$, with $\delta \in \{10^0, 10^{-1}, 10^{-2}, 10^{-3}\}$. The $x$-axis is the transformed variable of the imaginary part of $\omega'$ according to (18); 64 Legendre–Gauss–Radau points are used. (b) Zoom-in of the values of the integrand corresponding to different $\delta$.

## 5.3 Quadrature error

We use the Legendre–Gauss–Radau (LGR) quadrature rule with 8, 16, 32, and 64 quadrature points to evaluate the integral (19). These numerical integrals are denoted by $I_8$, $I_{16}$, $I_{32}$, and $I_{64}$, respectively. We estimate the errors in $I_8$, $I_{16}$ and $I_{32}$ by comparing them to $I_{64}$. For instance, the error in $I_8$ is estimated to be $|I_8 - I_{64}|$.

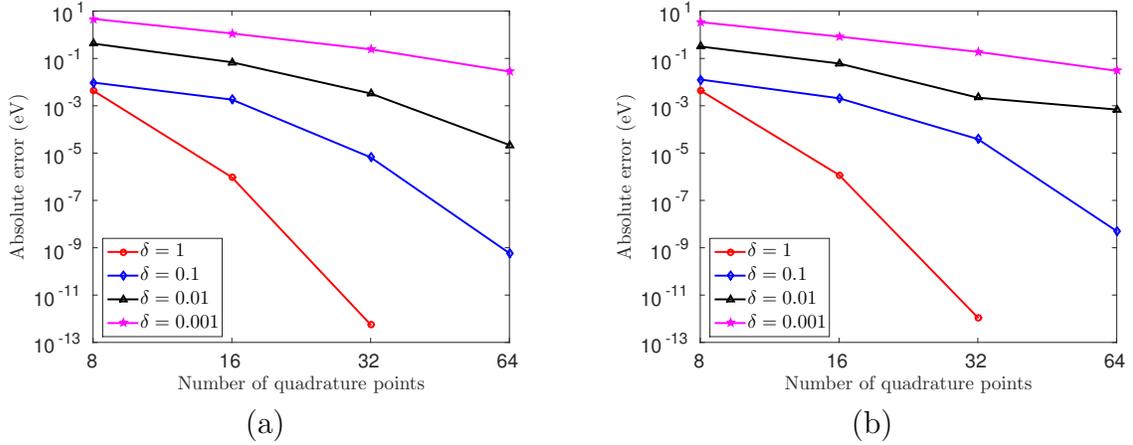

Figure 7: Absolute quadrature error resulting from applying the LGR rule to the evaluation of (19) along the vertical paths defined by (a) $\mathrm{Re}\,(\omega') = (1-\delta)\varepsilon_g/2$, (b) $\mathrm{Re}\,(\omega') = \delta_W - \delta\varepsilon_g/3$, with $\delta \in \{10^0, 10^{-1}, 10^{-2}, 10^{-3}\}$ for $\omega = (\varepsilon_{n_v} + \varepsilon_{n_v+1})/2$ in (a) and $\omega = 5\varepsilon_g/6 - \delta_W$ in (b) respectively.

Figure 7 (a) clearly shows that quadrature error increases when the integration path approaches a pole of $G_0$ and the number of quadrature points is fixed. Increasing the number of quadrature points can lead to significant improvement in accuracy when a good contour is chosen. Also, for a good contour the number of quadrature points required to obtain an accurate numerical integral can be as few as 8.

A similar pattern can be observed for contours that are progressively closer to a pole of $W_p$. Figure 7(b) shows the quadrature errors associated with applying the LGR rule with different quadrature points to contours in which the vertical paths are chosen along $\omega' = \delta_W - \delta\varepsilon_g/3 + i\zeta$, where $\zeta \in [0, +\infty)$, with $\delta \in \{10^0, 10^{-1}, 10^{-2}, 10^{-3}\}$, $\omega = 5\varepsilon_g/6 - \delta_W = -2.0641\,\mathrm{eV}$. The convergence of the numerical quadrature is similar to that shown in Figure 7(a).

## 5.4 Truncated SVD

Several low rank approximation schemes presented in Section 4 are based on finding a truncated SVD of $W_p$, $v\chi$ or $v\chi_0$. We would like the rank of the truncated SVD to be as small as possible without introducing too much error in the approximation. How much we can truncate depends on how fast the singular values of these operators decay.

In Figure 8, we plot all normalized singular values of $W_p$, $v\chi_0$, and $v\chi$, at $\omega' = 0$ and $\omega' = 2\mathrm{i}$. The normalization is performed by dividing all singular values by the largest one. We observe that the first dozen singular values decrease rapidly. Such rapid decrease makes it possible to construct a good low rank approximation of $W_p$, $v\chi_0$, and $v\chi$. We observe that the singular values of $W_p$ decreases the fastest, followed by those of $v\chi_0$, then by those of $v\chi$.

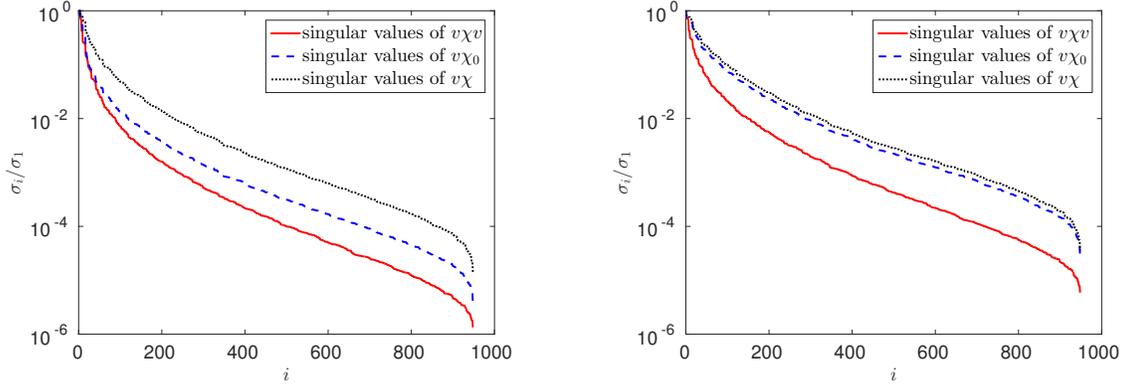

Figure 8: The normalized singular values of $v\chi v$, $v\chi_0$, and $v\chi$ at $\omega' = 0$ (left) and $\omega' = 2i$ (right).

## 5.5 Low rank approximation error

To examine the cumulative error introduced by both numerical integration and a low rank approximation to $W_p$, we plot in Figure 9 the real part of the integrand in (19) evaluated at 16 LGR quadrature points for $\omega = (\varepsilon_4 + \varepsilon_5)/2$. We used truncated SVD approximations to $W_p$ with different truncation levels, i.e., we approximate $W_p$ by $U_k S_k V_k^*$, where $S_k$ contains the $k$ largest singular values of $W_p$ and $U_k$ and $V_k$ contain the corresponding singular vectors. We can see from Figure 9 that the computed integrands are nearly indistinguishable when the rank of the SVD approximation is larger than 100.

In Figure 10, we plot the absolute error associated with the computed values of the numerical integral (19) obtained by using LGR with 64 quadrature points and different types of low rank approximations with different ranks. It is interesting to see that for a fixed rank, a low rank approximation obtained from a truncated SVD of $v\chi$ or $v\chi_0$ yields smaller error even though the singular values of $v\chi$ or $v\chi_0$ decreases less rapidly than those of $W_p$. It appears that truncating the bare Coulomb in $v\chi_0$ based on a small energy cutoff (see discussion in Section 4.2 produces much larger error.

We can also see from this figure that when the low rank approximation is con-structed from a truncated SVD of $v\chi_0$ that keeps only 50 largest singular values and the corresponding singular vectors, the absolute error in the computed correlation energy for the HOMO state is around 0.1 eV, which is generally acceptable. In this case, the number of operations required to construct the low rank approximation of $W_p$, which is the dominant cost of Algorithm 2 for each $\omega'$, is roughly $50n_v^3 n \approx 3 \times 10^6$. This is a lot smaller than the $n^2 n_v n_c \approx 3.3 \times 10^9$ operations required to construct $\chi_0$ in Algorithm 1.

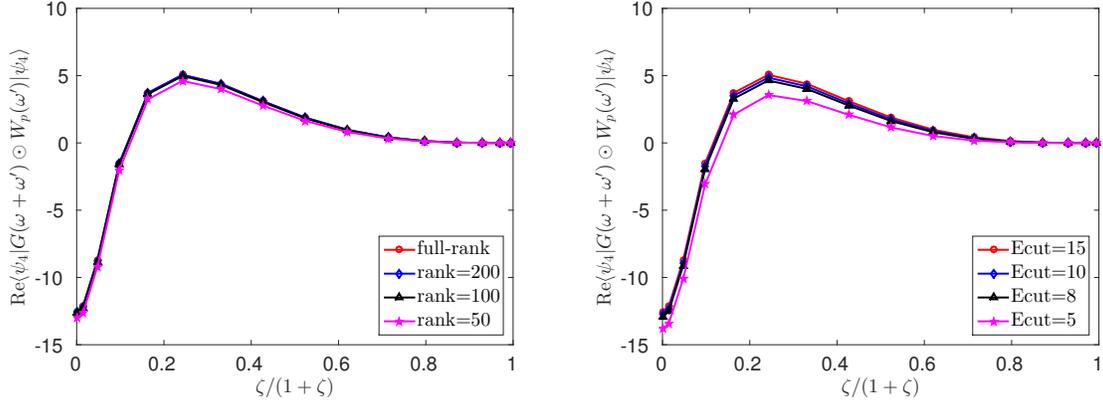

Figure 9: The integrand $\langle \psi_{n_v}^{\mathrm{KS}} | G_0(\omega + \omega') \odot W_p(\omega') | \psi_{n_v}^{\mathrm{KS}} \rangle$ evaluated on the imaginary axis (optimal contour) for $\omega = (\varepsilon_4 + \varepsilon_5)/2$. The contour integration is performed in the $\xi$ coordinate with 16 Legendre–Gauss–Radau (LGR) points. The $W_p$ is approximated by low rank approximations of the form (21) with different ranks (left), and low rank approximations based on Fourier truncation with different energy cutoff values (right).

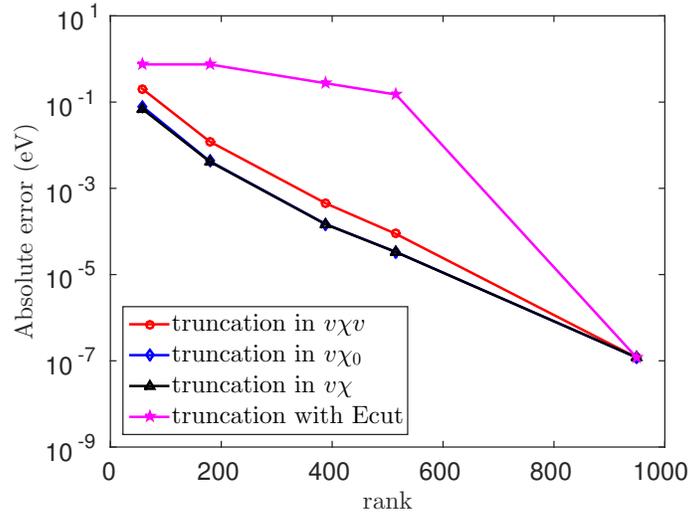

Figure 10: Absolute errors for the integral with low rank approximation.

Table 1: Computed results of $\Sigma_C$ (eV) with different truncation levels and different paths.

|  | full-rank | rank = 200 | rank = 100 | rank = 50 |
|---|---|---|---|---|
| $\delta = 1$ | $-1.4084$ | $-1.3999$ | $-1.3432$ | $-1.1538$ |
| $\delta = 10^{-3}$ | $-1.4084$ | $-1.4000$ | $-1.3427$ | $-1.1509$ |
| $\delta = -10^{-3}$ | $-1.4084$ | $-1.3999$ | $-1.3426$ | $-1.1509$ |

## 5.6 Low rank approximation on different contours

Finally, we perform a test to illustrate the effect of low rank approximation on different choices of contours. We still use the setting $\omega = (\varepsilon_4 + \varepsilon_5)/2$, and perform integration on the paths $\omega' = (1 - \delta)\varepsilon_g/2 + i\zeta$, where $\zeta \in (-\infty, +\infty)$, with $\delta \in \{10^0, 10^{-3}, -10^{-3}\}$. Truncated SVD on $W_p$ is applied in this test. We remark that for different integration paths, the same quadrature node in terms of $\xi$ (using the change of variable (18)) corresponds to different $\omega'$ with the same imaginary part.

Although we expect that choosing a contour close to a pole of $G_0$ may result in larger error in the integral, it turns out that there is very little difference among the numerical integrals performed along different contours, as we can see from Table 1, when 256 quadrature points are used in the numerical integration. We remark that for $\delta = -10^{-3}$, the contour actually encloses a pole of $G_0$. So the final result is adjusted by the residue term in (15).

Figure 11(a) illustrates the absolute error of the integrand in terms of the variable $\xi$ defined in (18). Despite that the integrands are quite different for different choices of contours (see Figure 6(a)), the errors are surprisingly close in the pointwise sense, except for a narrow region around $\xi = 0$, corresponding to points that are close to a pole of $G_0$. In Figure 11(b) we plot the integrand scaled by the Jacobian, $(1-\xi)^{-2} = (1+\zeta)^{-2}$, as well as the weight function of the LGR quadrature rule. The problematic region around $\xi = 0$ has a tiny weight and hence has limited contribution to the overall error of $\Sigma_C$. This briefly explains why the error introduced by low rank approximation can be insensitive to the choice of contour. However, the reason why the errors along different paths are similar is still not clear to us. We shall investigate this in detail in future work.

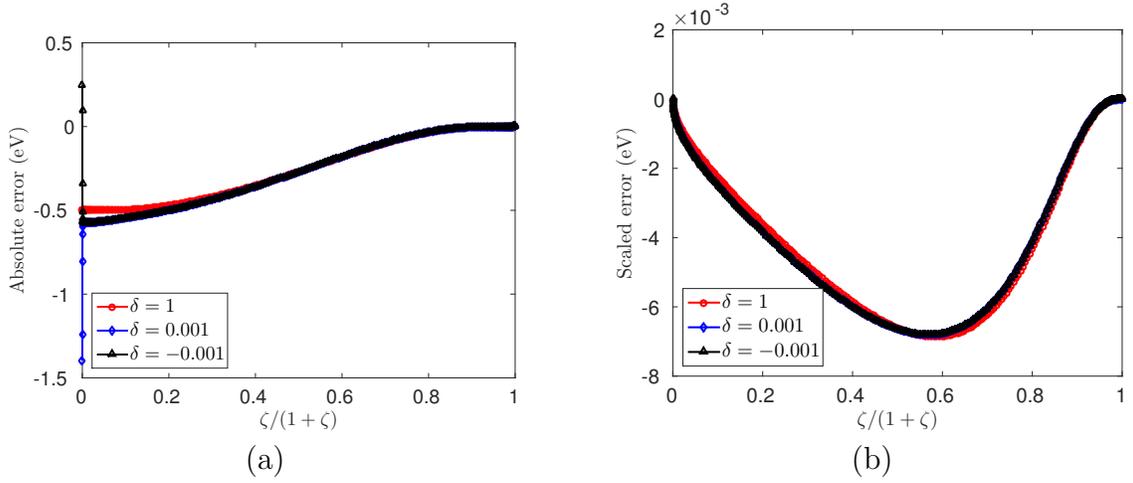

Figure 11: The absolute errors on 256 LGR quadrature nodes in (a) Re $\langle\psi_{n_v}^{\text{KS}}|G_0(\omega+\omega')\odot W_p(\omega')|\psi_{n_v}^{\text{KS}}\rangle$, (b) Re $\langle\psi_{n_v}^{\text{KS}}|G_0(\omega+\omega')\odot W_p(\omega')|\psi_{n_v}^{\text{KS}}\rangle w(\zeta)/(1+\zeta)^2$, introduced by low rank approximation, where $w(\zeta)$ is the weight function of the LGR quadrature rule. The $W_p$ is approximated by low rank approximations of the form (21) with rank = 50.

## 6 Summary

We presented a number of ways to construct a low rank approximation of the frequency dependent part of the screened Coulomb operator $W_p$, and examined the effect of low rank approximation on the accuracy of the computed correlation energy. We showed by a numerical example that we can keep the error introduced by a low rank approximation of $W_p$ under $0.1\,\text{eV}$ even when the rank of $W_p$ is an order of magnitude smaller than the dimension of $W_p$. We discussed how to choose an appropriate contour when using the contour deformation technique to evaluate the convolution between $G_0$ and $W_p$. We also showed how to apply the Legendre–Gauss–Radau (LGR) quadrature rule together with a special change of variable to perform numerical integration on a path parallel to the imaginary axis in the complex plane. The construction of a low rank approximation to $W_p$ and the subsequent computation that involves multiplying $G_0$ with a number of vectors requires solving a number of linear systems with a frequency shifted $H_{\text{KS}}$ as the coefficient matrix. We discussed how these linear systems can be solved by an iterative method. However, in our numerical experiments, we used the MATLAB SVD function and solved the linear system (36) using a spectral decomposition, in order to assess the error introduced exclusively by a low rank approximation of $W_p$. Solving these linear systems itera-tively is likely to introduce additional approximation errors that need to be carefully examine. We will investigate this type of error in future work. In addition, we will investigate efficient ways to solve these linear systems.

# Acknowledgments

This research was supported by the SciDAC Program on Excited State Phenomena in Energy Materials funded by the U. S. Department of Energy, Office of Basic Energy Sciences and of Advanced Scientific Computing Research, under Contract No. DE-AC02-05CH11231 at Lawrence Berkeley National Laboratory. L. L. acknowledges support from the Center for Applied Mathematics for Energy Research Applications (CAMERA) funded by U.S. Department of Energy, Office of Science, Advanced Scientific Computing Research and Basic Energy Sciences, and the Alfred P. Sloan fellowship. F. L. acknowledges support from the National Science Foundation of China under grant 11171232. Part of the research was conducted while she was a visiting scientist at LBNL.